\documentclass[12pt,leqno]{article}
\usepackage[a4paper, margin=25.4mm]{geometry}
\usepackage{amssymb}
\usepackage{amsmath}
\usepackage{amsthm}
\usepackage{stmaryrd}
\usepackage{mathtools,mathrsfs}
\usepackage[pagebackref]{hyperref}
\hypersetup{citecolor=purple, linkcolor=blue, colorlinks=true}
\usepackage{booktabs}
\usepackage{enumerate}
\usepackage{tikz}
\usepackage{url}
\usepackage[T2A,T1]{fontenc}
\usepackage{authblk}  

\theoremstyle{definition}

\renewcommand{\le}{\leqslant}
\renewcommand{\ge}{\geqslant}

\newcommand\scalemath[2]{\scalebox{#1}{\mbox{\ensuremath{\displaystyle #2}}}}

\definecolor{mygreen}{rgb}{0.55, 0.71, 0.0}
\definecolor{myyellow}{rgb}{1.0, 0.88, 0.21}
\definecolor{mywhite}{gray}{1.00}
\definecolor{mygrey}{gray}{0.92} 

\author[1]{S.P. Glasby}
\affil[1]{\small Center for the Mathematics of Symmetry and Computation,
  University of\newline Western Australia, Perth 6009, Australia\quad
  \href{mailto:Stephen.Glasby@uwa.edu.au}{Stephen.Glasby@uwa.edu.au}}

\title{Mathematics of the NYT daily word game Waffle}

\date{\today}

\begin{document}
\maketitle

\begin{abstract}
  This note investigates the combinatorics of permutations underlying
  the NYT daily word game Waffle. It helps to solve Waffle games and
  helps to understand why some games are easy to solve while others
  are very hard.  It shows that a perfect unscrambling must have
  precisely 11 orbits, with at least one of length 1, on the 21 Waffle
  squares. It also describes practical algorithms for solving Waffle
  games and creating new games with extreme properties.
\end{abstract}


\section{Introduction}

Word games have been very popular for a long time. In 1942 the New
York Times (NYT) started printed crosswords, initially to distract
readers from the bombing of Pearl Harbor, and later just for fun.
Since then many new daily/weekly word puzzles have emerged including
Mini crossword and Spelling bee in (2014). The two most popular NYT
daily word games are Wordle (2021/2), and Connections (2023).  These
games became popular at the start of the Covid pandemic, and are
anecdotally good for improving `brain function'.  This note explores
some mathematics (and computing) associated with the word game Waffle, which has
similar rules to Wordle. Both games were created by a software
engineers: Josh Wardle created Wordle in 2021, and James
Robinson~\cite{JR} created Waffle in 2022. For Wordle, the goal is to
find one five letter (English) word in six guesses, and for Waffle,
the goal is to find six five letter (English) words using 10 to 15
swaps. The words are scrambled and placed in a $5\times5$ grid with 4
`holes' which looks like a waffle. Over 1500 archived Waffle games are
found at~\cite{WA}. Figure~\ref{F:1}(a) is the scrambled archived~\cite{WA}
game \#2, and Figure~\ref{F:1}(b) is its unscrambled version.

\begin{figure}[!ht]
\begin{center}
\begin{tikzpicture}[scale=0.5] %
\draw [fill=mygrey] (0,0) rectangle (5,5);
\draw [fill=myyellow] (2,0) rectangle (3,2);
\draw [fill=myyellow] (3,0) rectangle (4,1);
\draw [fill=myyellow] (0,2) rectangle (1,3);
\draw [fill=myyellow] (4,2) rectangle (5,4);
\draw [fill=mygreen] (0,0) rectangle (1,1);
\draw [fill=mygreen] (4,0) rectangle (5,1);
\draw [fill=mygreen] (1,2) rectangle (3,3);
\draw [fill=mygreen] (0,3) rectangle (1,5);
\draw [fill=mygreen] (4,4) rectangle (5,5);
\draw [fill=mywhite] (1,1) rectangle (2,2);
\draw [fill=mywhite] (1,3) rectangle (2,4);
\draw [fill=mywhite] (3,1) rectangle (4,2);
\draw [fill=mywhite] (3,3) rectangle (4,4);
\draw [-] (0,0)-- (5,0);
\draw [-] (0,1)-- (5,1);
\draw [-] (0,2)-- (5,2);
\draw [-] (0,3)-- (5,3);
\draw [-] (0,4)-- (5,4);
\draw [-] (0,5)-- (5,5);
\draw [-] (0,0)-- (0,5);
\draw [-] (1,0)-- (1,5);
\draw [-] (2,0)-- (2,5);
\draw [-] (3,0)-- (3,5);
\draw [-] (4,0)-- (4,5);
\draw [-] (5,0)-- (5,5);
\node [align=center,scale=0.8,black] at (0.5,0.5) {F};
\node [align=center,scale=0.8,black] at (1.5,0.5) {F};
\node [align=center,scale=0.8,black] at (2.5,0.5) {A};
\node [align=center,scale=0.8,black] at (3.5,0.5) {R};
\node [align=center,scale=0.8,black] at (4.5,0.5) {E};
\node [align=center,scale=0.8,black] at (0.5,2.5) {I};
\node [align=center,scale=0.8,black] at (1.5,2.5) {N};
\node [align=center,scale=0.8,black] at (2.5,2.5) {D};
\node [align=center,scale=0.8,black] at (3.5,2.5) {E};
\node [align=center,scale=0.8,black] at (4.5,2.5) {E};
\node [align=center,scale=0.8,black] at (0.5,4.5) {S};
\node [align=center,scale=0.8,black] at (1.5,4.5) {C};
\node [align=center,scale=0.8,black] at (2.5,4.5) {G};
\node [align=center,scale=0.8,black] at (3.5,4.5) {O};
\node [align=center,scale=0.8,black] at (4.5,4.5) {L};
\node [align=center,scale=0.8,black] at (0.5,1.5) {R};
\node [align=center,scale=0.8,black] at (2.5,1.5) {I};
\node [align=center,scale=0.8,black] at (4.5,1.5) {U};
\node [align=center,scale=0.8,black] at (0.5,3.5) {N};
\node [align=center,scale=0.8,black] at (2.5,3.5) {N};
\node [align=center,scale=0.8,black] at (4.5,3.5) {D};
\end{tikzpicture}
\quad\quad\quad
\begin{tikzpicture}[scale=0.5] %
\draw [fill=mygrey] (0,0) rectangle (5,5);
\draw [fill=mygreen] (0,0) rectangle (5,1);
\draw [fill=mygreen] (0,2) rectangle (5,3);
\draw [fill=mygreen] (0,4) rectangle (5,5);
\draw [fill=mygreen] (0,0) rectangle (1,5);
\draw [fill=mygreen] (2,0) rectangle (3,5);
\draw [fill=mygreen] (4,0) rectangle (5,5);
\draw [-] (0,0)-- (5,0);
\draw [-] (0,1)-- (5,1);
\draw [-] (0,2)-- (5,2);
\draw [-] (0,3)-- (5,3);
\draw [-] (0,4)-- (5,4);
\draw [-] (0,5)-- (5,5);
\draw [-] (0,0)-- (0,5);
\draw [-] (1,0)-- (1,5);
\draw [-] (2,0)-- (2,5);
\draw [-] (3,0)-- (3,5);
\draw [-] (4,0)-- (4,5);
\draw [-] (5,0)-- (5,5);
\draw [fill=mywhite] (1,1) rectangle (2,2);
\draw [fill=mywhite] (1,3) rectangle (2,4);
\draw [fill=mywhite] (3,1) rectangle (4,2);
\draw [fill=mywhite] (3,3) rectangle (4,4);
\node [align=center,scale=0.8,black] at (0.5,0.5) {F};
\node [align=center,scale=0.8,black] at (1.5,0.5) {O};
\node [align=center,scale=0.8,black] at (2.5,0.5) {R};
\node [align=center,scale=0.8,black] at (3.5,0.5) {C};
\node [align=center,scale=0.8,black] at (4.5,0.5) {E};
\node [align=center,scale=0.8,black] at (0.5,2.5) {U};
\node [align=center,scale=0.8,black] at (1.5,2.5) {N};
\node [align=center,scale=0.8,black] at (2.5,2.5) {D};
\node [align=center,scale=0.8,black] at (3.5,2.5) {I};
\node [align=center,scale=0.8,black] at (4.5,2.5) {D};
\node [align=center,scale=0.8,black] at (0.5,4.5) {S};
\node [align=center,scale=0.8,black] at (1.5,4.5) {N};
\node [align=center,scale=0.8,black] at (2.5,4.5) {A};
\node [align=center,scale=0.8,black] at (3.5,4.5) {R};
\node [align=center,scale=0.8,black] at (4.5,4.5) {L};
\node [align=center,scale=0.8,black] at (0.5,1.5) {F};
\node [align=center,scale=0.8,black] at (2.5,1.5) {E};
\node [align=center,scale=0.8,black] at (4.5,1.5) {G};
\node [align=center,scale=0.8,black] at (0.5,3.5) {N};
\node [align=center,scale=0.8,black] at (2.5,3.5) {I};
\node [align=center,scale=0.8,black] at (4.5,3.5) {E};
\end{tikzpicture}
\quad\quad\quad
\begin{tikzpicture}[scale=0.5] %
\draw [fill=mygrey] (0,0) rectangle (5,5);
\draw [-] (0,0)-- (5,0);
\draw [-] (0,1)-- (5,1);
\draw [-] (0,2)-- (5,2);
\draw [-] (0,3)-- (5,3);
\draw [-] (0,4)-- (5,4);
\draw [-] (0,5)-- (5,5);
\draw [-] (0,0)-- (0,5);
\draw [-] (1,0)-- (1,5);
\draw [-] (2,0)-- (2,5);
\draw [-] (3,0)-- (3,5);
\draw [-] (4,0)-- (4,5);
\draw [-] (5,0)-- (5,5);
\draw [fill=mywhite] (1,1) rectangle (2,2);
\draw [fill=mywhite] (1,3) rectangle (2,4);
\draw [fill=mywhite] (3,1) rectangle (4,2);
\draw [fill=mywhite] (3,3) rectangle (4,4);
\node [align=center,scale=0.8,black] at (0.5,0.5) {17};
\node [align=center,scale=0.8,black] at (1.5,0.5) {18};
\node [align=center,scale=0.8,black] at (2.5,0.5) {19};
\node [align=center,scale=0.8,black] at (3.5,0.5) {20};
\node [align=center,scale=0.8,black] at (4.5,0.5) {21};
\node [align=center,scale=0.8,black] at (0.5,2.5) {9};
\node [align=center,scale=0.8,black] at (1.5,2.5) {10};
\node [align=center,scale=0.8,black] at (2.5,2.5) {11};
\node [align=center,scale=0.8,black] at (3.5,2.5) {12};
\node [align=center,scale=0.8,black] at (4.5,2.5) {13};
\node [align=center,scale=0.8,black] at (0.5,4.5) {1};
\node [align=center,scale=0.8,black] at (1.5,4.5) {2};
\node [align=center,scale=0.8,black] at (2.5,4.5) {3};
\node [align=center,scale=0.8,black] at (3.5,4.5) {4};
\node [align=center,scale=0.8,black] at (4.5,4.5) {5};
\node [align=center,scale=0.8,black] at (0.5,1.5) {14};
\node [align=center,scale=0.8,black] at (2.5,1.5) {15};
\node [align=center,scale=0.8,black] at (4.5,1.5) {16};
\node [align=center,scale=0.8,black] at (0.5,3.5) {6};
\node [align=center,scale=0.8,black] at (2.5,3.5) {7};
\node [align=center,scale=0.8,black] at (4.5,3.5) {8};
\end{tikzpicture}
\quad\quad\quad
\begin{tikzpicture}[scale=0.5] 
\draw [fill=mygrey] (0,0) rectangle (5,5);
\draw [-] (0,0)-- (5,0);
\draw [-] (0,1)-- (5,1);
\draw [-] (0,2)-- (5,2);
\draw [-] (0,3)-- (5,3);
\draw [-] (0,4)-- (5,4);
\draw [-] (0,5)-- (5,5);
\draw [-] (0,0)-- (0,5);
\draw [-] (1,0)-- (1,5);
\draw [-] (2,0)-- (2,5);
\draw [-] (3,0)-- (3,5);
\draw [-] (4,0)-- (4,5);
\draw [-] (5,0)-- (5,5);
\draw [fill=myyellow] (0,0) rectangle (1,1);
\draw [fill=myyellow] (2,0) rectangle (3,1);
\draw [fill=myyellow] (4,0) rectangle (5,1);
\draw [fill=myyellow] (0,2) rectangle (1,3);
\draw [fill=myyellow] (2,2) rectangle (3,3);
\draw [fill=myyellow] (4,2) rectangle (5,3);
\draw [fill=myyellow] (0,4) rectangle (1,5);
\draw [fill=myyellow] (2,4) rectangle (3,5);
\draw [fill=myyellow] (4,4) rectangle (5,5);
\draw [fill=mywhite] (1,1) rectangle (2,2);
\draw [fill=mywhite] (1,3) rectangle (2,4);
\draw [fill=mywhite] (3,1) rectangle (4,2);
\draw [fill=mywhite] (3,3) rectangle (4,4);
\node [align=center,scale=0.8,black] at (0.5,0.5) {O};
\node [align=center,scale=0.8,black] at (1.5,0.5) {E};
\node [align=center,scale=0.8,black] at (2.5,0.5) {O};
\node [align=center,scale=0.8,black] at (3.5,0.5) {E};
\node [align=center,scale=0.8,black] at (4.5,0.5) {O};
\node [align=center,scale=0.8,black] at (0.5,2.5) {O};
\node [align=center,scale=0.8,black] at (1.5,2.5) {E};
\node [align=center,scale=0.8,black] at (2.5,2.5) {O};
\node [align=center,scale=0.8,black] at (3.5,2.5) {E};
\node [align=center,scale=0.8,black] at (4.5,2.5) {O};
\node [align=center,scale=0.8,black] at (0.5,4.5) {O};
\node [align=center,scale=0.8,black] at (1.5,4.5) {E};
\node [align=center,scale=0.8,black] at (2.5,4.5) {O};
\node [align=center,scale=0.8,black] at (3.5,4.5) {E};
\node [align=center,scale=0.8,black] at (4.5,4.5) {O};
\node [align=center,scale=0.8,black] at (0.5,1.5) {E};
\node [align=center,scale=0.8,black] at (2.5,1.5) {E};
\node [align=center,scale=0.8,black] at (4.5,1.5) {E};
\node [align=center,scale=0.8,black] at (0.5,3.5) {E};
\node [align=center,scale=0.8,black] at (2.5,3.5) {E};
\node [align=center,scale=0.8,black] at (4.5,3.5) {E};
\end{tikzpicture}
\vskip-2mm
\caption{(a) scrambled game \#2 (b) unscrambled (c) numbered squares (d) parities}
\label{F:1}
\end{center}
\end{figure}

The coloring of a square, like in Wordle, gives
information about whether a given letter is in the correct place (green), or is
incorrect but in the word (yellow), or is not in the word (gray).
We discuss coloring in Section~\ref{S:color}.
A successful Waffle game finds the (unique) solution in at most 15 swaps.
It turns out that every game can be uniquely solved with 10 swaps and no fewer,
and a perfect ten-swap-game earns a maximum of 5 points. I will assume that
the reader is curious about the mathematics to achieve a perfect score.
Section~\ref{Perfect} describes the types of `perfect unscrambling' that
arise. Sections~\ref{S:color} and~\ref{S:hard} consider unique solutions, and
what makes a Waffle game hard, respectively.

\section{Perfect unscramblings and Cayley's lemma}\label{Perfect}

Three steps are needed to obtain a perfect score (5 out of 5 with 10 swaps) in Waffle:
\begin{enumerate}[Step 1.]
  \setlength{\topsep}{2pt}\setlength{\itemsep}{0pt}\setlength{\parskip}{0pt}
  \item find the six 5-letter words, 
  \item find a perfect unscrambling, and
  \item find 10 swaps of the squares which give the perfect unscrambling.
\end{enumerate}
\kern-7pt
Step~3 is easy to solve using ideas below. This section shows that
`many' repeat letters can make Step~2 hard, and no repeat letters can make
Step~1 hard. Allowing mistakes makes Steps~1,2 both easier, but we seek
a perfect score of~5, with no mistakes and~10~swaps.

To obtain a perfect score in Waffle, you need a \emph{perfect unscrambling}.
This is a permutation of the $5^2-4=21$ colored squares obtained using ten
swaps (or transpositions) that will unscramble the puzzle to give
the six 5-letter words. Imagine that we are permuting (or rearranging) the
21 \emph{numbered squares} in Figure~\ref{F:1}(c) rather than the \emph{letters}
which may contain repeats. The symmetric group $S_{21}$
contains $21!$ permutations of 21 distinct objects.
But the perfect unscramblings are
permutations that can be obtained using 10 swaps \emph{and no fewer} (recall that no game can be solved in fewer than 10 swaps). What do these look like?
We will prove that perfect unscramblings are atypical permutations that
contain more than the expected number of cycles:
they have precisely 11 disjoint cycles whereas the average number of
disjoint cycles for a uniformly random $g\in S_{21}$ is $\approx3.65$.
This average was computed in~\cite{G}, by considering the 792 conjugacy
classes of $S_{21}$ and their sizes.

\begin{figure}[!ht]
  \begin{center}
  \begin{tikzpicture}
   \coordinate (node0) at (0:1/2);
   \draw [fill] (node0) circle [radius=0.05];
   \node [below] at (node0) {$\scalemath{0.8}{x_1}$};
   \begin{scope}[xshift=2cm]
    \foreach \i in {0, 1} {
       \coordinate (node\i) at (360/2*\i:1/2);
       \draw [fill] (node\i) circle [radius=0.05];
     }
     \draw (node0) -- (node1);
     \node [below] at (node0) {$\scalemath{0.8}{x_2}$};
     \node [below] at (node1) {$\scalemath{0.8}{x_3}$};
   \end{scope}
   
   \begin{scope}[xshift=4cm]
       \foreach \i in {0, 1, 2} {
           \coordinate (node\i) at (360/3*\i:1/2);
           \draw [fill] (node\i) circle [radius=0.05];
       }
       \draw (node0) -- (node1);
       \draw (node1) -- (node2);
       \draw (node2) -- (node0);
       \node [right] at (node0) {$\scalemath{0.8}{x_4}$};
       \node [left] at (node1) {$\scalemath{0.8}{x_5}$};
       \node [below] at (node2) {$\scalemath{0.8}{x_6}$};
   \end{scope}
   
   \begin{scope}[xshift=6.5cm]
       \foreach \i in {0, 1, 2, 3} {
           \coordinate (node\i) at (360/4*\i:1/2);
           \draw [fill] (node\i) circle [radius=0.05];
       }
       \draw (node0) -- (node1);
       \draw (node1) -- (node2);
       \draw (node2) -- (node3);
       \draw (node3) -- (node0);
       \node [right] at (node0) {$\scalemath{0.8}{x_7}$};
       \node [left] at (node1) {$\scalemath{0.8}{x_8}$};
       \node [left] at (node2) {$\scalemath{0.8}{x_9}$};
       \node [right] at (node3) {$\scalemath{0.8}{x_{10}}$};
   \end{scope}
   
   \begin{scope}[xshift=9cm]
       \foreach \i in {0, 1, 2, 3, 4} {
      \coordinate (node\i) at (360/5*\i:1/2);
           \draw [fill] (node\i) circle [radius=0.05];
       }
       \draw (node0) -- (node1);
       \draw (node1) -- (node2);
       \draw (node2) -- (node3);
       \draw (node3) -- (node4);
       \draw (node4) -- (node0);
       \node [right] at (node0) {$\scalemath{0.8}{x_{11}}$};
       \node [above] at (node1) {$\scalemath{0.8}{x_{12}}$};
       \node [left] at (node2) {$\scalemath{0.8}{x_{13}}$};
       \node [left] at (node3) {$\scalemath{0.8}{x_{14}}$};
       \node [below] at (node4) {$\scalemath{0.8}{x_{15}}$};
   \end{scope}
   
   \begin{scope}[xshift=11.5cm, yshift=0cm]
      \foreach \i in {0, 1, 2, 3, 4, 5} {
      \coordinate (node\i) at (360/6*\i:1/2);
           \draw [fill] (node\i) circle [radius=0.05];
       }
       \draw (node0) -- (node1);
       \draw (node1) -- (node2);
       \draw (node2) -- (node3);
       \draw (node3) -- (node4);
       \draw (node4) -- (node5);
       \draw (node5) -- (node0);
       \node [right] at (node0) {$\scalemath{0.8}{x_{16}}$};
       \node [above] at (node1) {\ $\scalemath{0.8}{x_{17}}$};
       \node [above] at (node2) {$\scalemath{0.8}{x_{18}}$\ };
       \node [left] at (node3) {$\scalemath{0.8}{x_{19}}$};
       \node [below] at (node4) {$\scalemath{0.8}{x_{20}}$\ };
       \node [below] at (node5) {\ $\scalemath{0.8}{x_{21}}$};
   \end{scope}
  \end{tikzpicture}
  \end{center}
  \vskip-8mm
  \caption{$\scalemath{0.9}{g=(x_1)(x_2,x_3)(x_4,x_5,x_6)(x_7,x_8,x_9,x_{10})(x_{11},x_{12},x_{13},x_{14},x_{15})(x_{16},x_{17},x_{18},x_{19},x_{20},x_{21})}$}
  \label{F:1.5}
\end{figure}
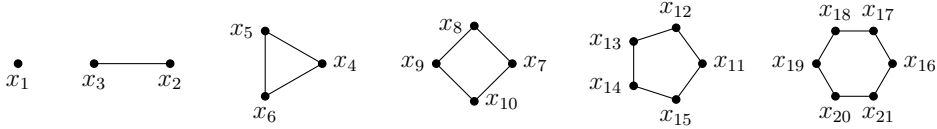

The \emph{symmetric group} $S_n$ comprises the $n!$ permutations of a set
$X=\{x_1,\dots,x_n\}$ with $n$ elements. Every permutation $g\in S_n$ can be
`factored' as a product of disjoint cycles, just as
whole numbers can be factored as a product of primes, see Figure~\ref{F:1.5}.
Let $c(g)$ be the \emph{number of disjoint cycles} of $g$, and let $s(g)$ be the
\emph{minimal number of swaps} needed form $g$.
For example, if $g\in S_{21}$ is the permutation
\setcounter{MaxMatrixCols}{21}
\[
\scalemath{0.87}{
\begin{pmatrix}
  x_1&x_2&x_3&x_4&x_5&x_6&x_7&x_8&x_9&x_{10}&x_{11}&x_{12}&x_{13}&x_{14}&x_{15}&x_{16}&x_{17}&x_{18}&x_{19}&x_{20}&x_{21}\\
  x_1&x_3&x_2&x_5&x_6&x_4&x_8&x_9&x_{10}&x_7&x_{12}&x_{13}&x_{14}&x_{15}&x_{11}&x_{17}&x_{18}&x_{19}&x_{20}&x_{21}&x_{16}
\end{pmatrix},}
\]
illustrated in Figure~\ref{F:1.5} where the disjoint cycles of $g$ turn
\emph{counterclockwise}, e.g. $g$ maps $x_4$ to $x_5$, $x_5$ to $x_6$
and $x_6$ back to $x_4$. 
Thus $g$ is a product of six cycles (of lengths 1,\,\dots,\,6) so $c(g)=6$.
Since a cycle of length $k$ can be formed from $k-1$ swaps and no fewer it follows that $s(g)=0+1+2+3+4+5=15$. (Can you formally prove that $s(g)=15$?)

A beautiful lemma of Arthur Cayley says that
$s(g)=n-c(g)$ for each $g\in S_n$.  Arthur Cayley (1821--1895) wrote
over 300 papers when he was a lawyer before becoming a Sadleirian
professor in mathematics at Cambridge in 1863, and he seemed to excel at
most everything. A result due to Gon\v{c}arov in 1944, see~\cite[p.\,176,
  Eq.\,(2.2)]{ET2}, says $c(g)\sim\log(n)$ as $n\to\infty$.
It is not hard to prove that the average value of $c(g)$ over all $g\in S_{21}$
is the $n$th harmonic number $H_n=\sum_{i=1}^n1/i$ and
$\log(n)+\gamma<H_n<\log(n)+\gamma+1/(2n)$ where $\gamma=0.57721\cdots$
is the Euler-Mascheroni constant $\lim_{n\to\infty}(H_n-\log n)$. When $n=21$ this average
is $\approx 3.65$, so that an average permutation is far from
a perfect unscrambling.

Cayley's lemma is related to graphs.
Let $G=S_n$, and let $S$ be the set of $\binom{n}{2}$ transpositions.
The Cayley graph $\Gamma_n=\textup{Cay}(G,S)$ has vertex set $V=G$ and
edge set
\begin{tikzpicture}[yshift=0mm]
  \draw [fill] (0,0) circle [radius=0.05];
  \draw [fill] (0.7,0) circle [radius=0.05];
  \draw (0,0)--(0.7,0);
  \node [above] at (0,0) {$\scalemath{0.7}{g}$};
  \node [above] at (0.7,0) {$\scalemath{0.7}{gs}$};
\end{tikzpicture}
where $g\in G$ and $s\in S$.
The distance $d(g,h)$ between
$g,h\in G$ is the length of the smallest path from $g$ to $h$ in
$\Gamma_n$. That is, $d(g,h)=s(g^{-1}h)$ is the fewest number of swaps
needed to change $g$ to $h$. We sketch a proof of Cayley's lemma:
If the disjoint cycle
decomposition of $g$ is $g_1\cdots g_c$ where $g_i\in S_{n_i}$ is a cycle of
length $n_i$ for each $i$, then
\[
  s(g)=\sum_{i=1}^c s(g_i)=\sum_{i=1}^c (n_i-1)=n-c\qquad
  \textup{where $n=\sum_{i=1}^cn_i$,}
\]
and the $n_i$-cycle $g_i=(x_1,\dots,x_{n_i})\in S_{n_i}$ has $s(g)=n_i-1$.
Note that
$s(g_i)\le n_i-1$ because $g_i$ is obtained from the $n_i-1$ swaps $(x_{n_i},x_{n_i-1})$ then
$\cdots$ then $(x_3,x_2)$ then $(x_2,x_1)$. And if fewer than $n_i-1$ swaps are
used, then fewer than $n_i$ points are moved.
  
We stress that $c(g)$ counts \emph{all} cycles for $g\in S_n$, including 1-cycles (a.k.a. fixed points).
Arguing as above, an `average' $g\in S_n$ is a product
of about $n-\log(n)$ swaps and no fewer. Interestingly, the largest
distance $d(g,h)$ between $g,h\in S_n$ is $n-1$ precisely
when $g^{-1}h$ is an $n$-cycle. (Clearly $s(g^{-1}h)$ is a maximum
when $c(g^{-1}h)=n-s(g^{-1}h)$ is a minimum.)
Hence the Cayley graph $\Gamma_n$ has diameter is $n-1$,
degree $\binom{n}{2}$ and order~$n!$.
  
It is important for us that the following are equivalent:
\begin{enumerate}[{\rm 1.}]
  \setlength{\topsep}{2pt}\setlength{\itemsep}{0pt}\setlength{\parskip}{0pt}
  \item~$g\in S_{21}$ is a perfect unscrambling for a Waffle game,
  \item~$g$ can be formed from $s(g)=10$ (not necessarily disjoint) swaps and no fewer, and
  \item~$g$ has disjoint cycle decomposition
    $g_1g_2\cdots g_{11}$ with $c(g)=21-s(g)=11$ cycles.
\end{enumerate}
Note that ten swaps will give a permutation
moving at most $2\times10=20$ squares, so each perfect unscrambling
$g\in S_{21}$ has at least 1 fixed (green) square, and at most
10 fixed (green) squares (because $c(g)=11$). The archived
Waffle games~\cite{WA} have green squares at the corners and center
(namely at positions $1, 5, 17, 21$ and $11$). This is not a \emph{requirement}
for Waffle, but having 5 to 8 green squares makes the game not too
hard or easy. We investigate in Section~\ref{S:color} what happens if there
is precisely 1 green square, the minimum number of fixed squares.
By the way, averaging the number of fixed points over all permutations
$g\in S_n$ gives 1. (Note that $(n-1)!$
elements of $S_n$ fix a given point, so the probability that a given
point is fixed is $1/n$. This gives an average of 1 fixed point.) 

We now measure how hard it can be to find a perfect unscrambling.
Let $g$ be the perfect unscrambling of $n=21$ squares for the archived Waffle
game \#2 in Figure~\ref{F:1}(a). Of its $n$ squares, precisely 7 are green,
and hence are fixed by $g$. Let $g'$ be the permutation induced by $g$ on the
remaining $n'=21-7=14$ non-green
squares. Each non-green letter appears once except ${\rm I, E, R}$ appear twice.
To distinguish repeats let ${\rm I}_1$ be in position 9, ${\rm I}_2$ in
position 15, etc. The six words in Figure~\ref{F:1}(b) show that $g'$
permutes the $n'=14$ non-green squares as follows:
\setcounter{MaxMatrixCols}{15}
\[
g'=
\begin{pmatrix}
  {\rm C}&{\rm G}&{\rm O}&{\rm N}&{\rm D}&{\rm I}_1&{\rm E}_1&{\rm E}_2&{\rm R}_1&{\rm I}_2&{\rm U}&{\rm F}&{\rm A}&{\rm R}_2\\
  {\rm N}&{\rm A}&{\rm R}_*&{\rm I}_*&{\rm E}_*&{\rm U}&{\rm I}_*&{\rm D}&{\rm F}&{\rm E}_*&{\rm G}&{\rm O}&{\rm R}_*&{\rm C}
  \end{pmatrix}
\]
where $*$ denotes a subscript of 1 or 2. As $g'$ must be
surjective (onto) the 14 symbols on the top line must appear on the bottom.
Therefore there are $2^3=8$ possible permutations for $g'$, and hence for $g$.
Applying Cayley's lemma to $g$ and $g'$ gives $s(g)=n-c(g)$ and
$s(g')=n'-c(g')$. However, $n=7+n'$ and $c(g)=7+c(g')$, so that $s(g)=s(g')=10$.
As $n=21$, we have $c(g)=11$ and $c(g')=4$.
Of the 8 possibilities for $g'$ only the following permutation has 4 cycles
\[
g'=
\begin{pmatrix}
  {\rm C}&{\rm G}&{\rm O}&{\rm N}&{\rm D}&{\rm I}_1&{\rm E}_1&{\rm E}_2&{\rm R}_1&{\rm I}_2&{\rm U}&{\rm F}&{\rm A}&{\rm R}_2\\
  {\rm N}&{\rm A}&{\rm R}_1&{\rm I}_1&{\rm E}_2&{\rm U}&{\rm I}_2&{\rm D}&{\rm F}&{\rm E}_1&{\rm G}&{\rm O}&{\rm R}_2&{\rm C}
  \end{pmatrix};
\]
the other 7 have fewer cycles. Indeed $g'=({\rm C}{\rm N}{\rm
  I}_1{\rm U}{\rm G}{\rm A}{\rm R}_2)({\rm O}{\rm R}_1{\rm F})({\rm
  D}{\rm E}_2)({\rm E}_1{\rm I}_2)$ is its disjoint cycle
decomposition.  This shows that we must determine more than the six
words in Figure~\ref{F:1}(b): we must also find the perfect
unscrambling $g$ in $S_{21}$ amongst the possibilities arising from
repeated letters, and then find a sequence of 10 swaps
giving~$g$. Similarly, in archived Waffle~\cite{WA} game \#7, there
are 5 green squares, and the permutation $g'$ on the remaining 16
squares has distinct letters except for 5 Rs and 2 Ns, so
$5!\cdot2!=240$ choices for $g'$. Only one has $c(g')=c(g)-5=6$ cycles, all
others have fewer cycles, so in this case finding a perfect unscrambling for
this puzzle by guessing has a 1/240 chance of success, even if you know the
correct placement of
all the letters (minus their distinguishing subscripts). In personal
communication, Carl Mueller~\cite{CM} 
described to me some helpful heuristics (without a formal proof) for finding
a perfect unscrambling (Step~2) very efficiently. The reason these heuristics
are likely to work is described below.

In general, if a perfect unscrambling $g\in S_{21}$ has $N_g$ fixed (green) squares,
then the permutation $g'$ of the remaining $n'=21-N_g$ squares has no fixed
squares. As $c(g)=11$, we have $c(g')=11-N_g$. Many archived Waffle games~\cite{WA} have
$N_g=6$, in this case $g'\in S_{15}$ has $c(g')=5$ whereas an `average' element of
$S_{15}$ has only $\approx3.32$ cycles by~\cite{G}. This reasoning suggests why
so many possible choices for $g'$ have fewer than $11-N_g$~cycles.

\section{Coloring and uniqueness}\label{S:color}

So far we have only considered green squares (which are in the right place).
To discuss yellow and gray squares we introduce the \emph{parity} of a square.
A square in row $i$ and column $j$ is called \emph{even} if the number
$10(i-1)+j$ is even, and \emph{odd} otherwise, see Figure~\ref{F:1}(d).
Hence the $3^2=9$ odd squares have $i,j\in\{1,3,5\}$
and they lie at the intersection of a horizontal and a vertical word. (Vertical
words are read from top to bottom with increasing column numbers.)
By contrast, the 12 even squares lie in only one word. The $2^2=4$ holes
have $i,j\in\{2,4\}$, so an even square is in a horizontal word
if $j\in\{2,4\}$, and in a vertical word if $i\in\{2,4\}$. 
An even square is colored gray if it does not belong to the unique word
incident with this square, and is yellow otherwise. An odd square is colored
gray if it does not belong to the two words incident with this square, and
is yellow otherwise. Waffle's gray coloring rules
\url{https://wafflegame.net/daily} are more complicated
for repeated letters, but this nuance will not concern us.

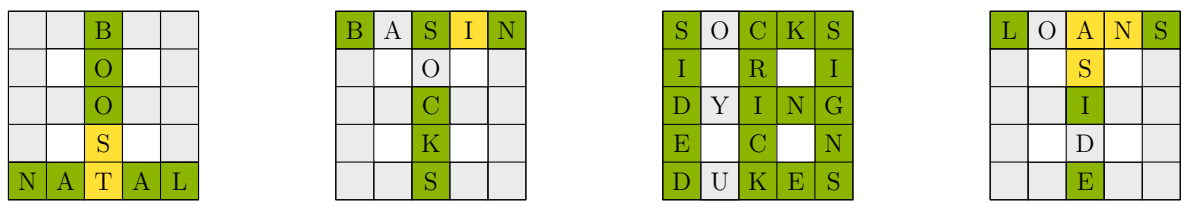
\begin{figure}[!ht]
\begin{center}
\begin{tikzpicture}[scale=0.5] %
\draw [fill=mygreen] (0,0) rectangle (2,1);
\draw [fill=mygreen] (3,0) rectangle (5,1);
\draw [fill=mygreen] (2,2) rectangle (3,5);
\draw [fill=mygrey] (0,1) rectangle (1,5);
\draw [fill=mygrey] (4,1) rectangle (5,5);
\draw [fill=mygrey] (1,2) rectangle (2,3);
\draw [fill=mygrey] (3,2) rectangle (4,3);
\draw [fill=mygrey] (1,4) rectangle (2,5);
\draw [fill=mygrey] (3,4) rectangle (4,5);
\draw [fill=myyellow] (2,0) rectangle (3,2);
\draw [-] (0,0)-- (5,0);
\draw [-] (0,1)-- (5,1);
\draw [-] (0,2)-- (5,2);
\draw [-] (0,3)-- (5,3);
\draw [-] (0,4)-- (5,4);
\draw [-] (0,5)-- (5,5);
\draw [-] (0,0)-- (0,5);
\draw [-] (1,0)-- (1,5);
\draw [-] (2,0)-- (2,5);
\draw [-] (3,0)-- (3,5);
\draw [-] (4,0)-- (4,5);
\draw [-] (5,0)-- (5,5);
\draw [fill=mywhite] (1,1) rectangle (2,2);
\draw [fill=mywhite] (1,3) rectangle (2,4);
\draw [fill=mywhite] (3,1) rectangle (4,2);
\draw [fill=mywhite] (3,3) rectangle (4,4);
\node [align=center,scale=0.8,black] at (0.5,0.5) {N};
\node [align=center,scale=0.8,black] at (1.5,0.5) {A};
\node [align=center,scale=0.8,black] at (2.5,0.5) {T};
\node [align=center,scale=0.8,black] at (3.5,0.5) {A};
\node [align=center,scale=0.8,black] at (4.5,0.5) {L};
\node [align=center,scale=0.8,black] at (0.5,2.5) {};
\node [align=center,scale=0.8,black] at (1.5,2.5) {};
\node [align=center,scale=0.8,black] at (2.5,2.5) {O};
\node [align=center,scale=0.8,black] at (3.5,2.5) {};
\node [align=center,scale=0.8,black] at (4.5,2.5) {};
\node [align=center,scale=0.8,black] at (0.5,4.5) {};
\node [align=center,scale=0.8,black] at (1.5,4.5) {};
\node [align=center,scale=0.8,black] at (2.5,4.5) {B};
\node [align=center,scale=0.8,black] at (3.5,4.5) {};
\node [align=center,scale=0.8,black] at (4.5,4.5) {};
\node [align=center,scale=0.8,black] at (0.5,1.5) {};
\node [align=center,scale=0.8,black] at (2.5,1.5) {S};
\node [align=center,scale=0.8,black] at (4.5,1.5) {};
\node [align=center,scale=0.8,black] at (0.5,3.5) {};
\node [align=center,scale=0.8,black] at (2.5,3.5) {O};
\node [align=center,scale=0.8,black] at (4.5,3.5) {};
\end{tikzpicture}
\qquad\qquad
\begin{tikzpicture}[scale=0.5] %
\draw [fill=mygreen] (0,4) rectangle (5,5);
\draw [fill=mygreen] (2,0) rectangle (3,3);
\draw [fill=mygrey] (1,4) rectangle (2,5); 
\draw [fill=myyellow] (3,4) rectangle (4,5);
\draw [fill=mygrey] (2,3) rectangle (3,4); 
\draw [fill=mygrey] (0,0) rectangle (1,4);
\draw [fill=mygrey] (4,0) rectangle (5,4);
\draw [fill=mygrey] (1,0) rectangle (2,1);
\draw [fill=mygrey] (3,0) rectangle (4,1);
\draw [fill=mygrey] (1,2) rectangle (2,3);
\draw [fill=mygrey] (3,2) rectangle (4,3);
\draw [-] (0,0)-- (5,0);
\draw [-] (0,1)-- (5,1);
\draw [-] (0,2)-- (5,2);
\draw [-] (0,3)-- (5,3);
\draw [-] (0,4)-- (5,4);
\draw [-] (0,5)-- (5,5);
\draw [-] (0,0)-- (0,5);
\draw [-] (1,0)-- (1,5);
\draw [-] (2,0)-- (2,5);
\draw [-] (3,0)-- (3,5);
\draw [-] (4,0)-- (4,5);
\draw [-] (5,0)-- (5,5);
\draw [fill=mywhite] (1,1) rectangle (2,2);
\draw [fill=mywhite] (1,3) rectangle (2,4);
\draw [fill=mywhite] (3,1) rectangle (4,2);
\draw [fill=mywhite] (3,3) rectangle (4,4);
\node [align=center,scale=0.8,black] at (0.5,0.5) {};
\node [align=center,scale=0.8,black] at (1.5,0.5) {};
\node [align=center,scale=0.8,black] at (2.5,0.5) {S};
\node [align=center,scale=0.8,black] at (3.5,0.5) {};
\node [align=center,scale=0.8,black] at (4.5,0.5) {};
\node [align=center,scale=0.8,black] at (0.5,2.5) {};
\node [align=center,scale=0.8,black] at (1.5,2.5) {};
\node [align=center,scale=0.8,black] at (2.5,2.5) {C};
\node [align=center,scale=0.8,black] at (3.5,2.5) {};
\node [align=center,scale=0.8,black] at (4.5,2.5) {};
\node [align=center,scale=0.8,black] at (0.5,4.5) {B};
\node [align=center,scale=0.8,black] at (1.5,4.5) {A};
\node [align=center,scale=0.8,black] at (2.5,4.5) {S};
\node [align=center,scale=0.8,black] at (3.5,4.5) {I};
\node [align=center,scale=0.8,black] at (4.5,4.5) {N};
\node [align=center,scale=0.8,black] at (0.5,1.5) {};
\node [align=center,scale=0.8,black] at (2.5,1.5) {K};
\node [align=center,scale=0.8,black] at (4.5,1.5) {};
\node [align=center,scale=0.8,black] at (0.5,3.5) {};
\node [align=center,scale=0.8,black] at (2.5,3.5) {O};
\node [align=center,scale=0.8,black] at (4.5,3.5) {};
\end{tikzpicture}
\qquad\qquad
\begin{tikzpicture}[scale=0.5] %
\draw [fill=mygreen] (0,0) rectangle (5,5);
\draw [fill=mygrey] (1,4) rectangle (2,5);
\draw [fill=mygrey] (1,2) rectangle (2,3);
\draw [fill=mygrey] (1,0) rectangle (2,1);
\draw [-] (0,0)-- (5,0);
\draw [-] (0,1)-- (5,1);
\draw [-] (0,2)-- (5,2);
\draw [-] (0,3)-- (5,3);
\draw [-] (0,4)-- (5,4);
\draw [-] (0,5)-- (5,5);
\draw [-] (0,0)-- (0,5);
\draw [-] (1,0)-- (1,5);
\draw [-] (2,0)-- (2,5);
\draw [-] (3,0)-- (3,5);
\draw [-] (4,0)-- (4,5);
\draw [-] (5,0)-- (5,5);
\draw [fill=mywhite] (1,1) rectangle (2,2);
\draw [fill=mywhite] (1,3) rectangle (2,4);
\draw [fill=mywhite] (3,1) rectangle (4,2);
\draw [fill=mywhite] (3,3) rectangle (4,4);
\node [align=center,scale=0.8,black] at (0.5,0.5) {D};
\node [align=center,scale=0.8,black] at (1.5,0.5) {U};
\node [align=center,scale=0.8,black] at (2.5,0.5) {K};
\node [align=center,scale=0.8,black] at (3.5,0.5) {E};
\node [align=center,scale=0.8,black] at (4.5,0.5) {S};
\node [align=center,scale=0.8,black] at (0.5,2.5) {D};
\node [align=center,scale=0.8,black] at (1.5,2.5) {Y};
\node [align=center,scale=0.8,black] at (2.5,2.5) {I};
\node [align=center,scale=0.8,black] at (3.5,2.5) {N};
\node [align=center,scale=0.8,black] at (4.5,2.5) {G};
\node [align=center,scale=0.8,black] at (0.5,4.5) {S};
\node [align=center,scale=0.8,black] at (1.5,4.5) {O};
\node [align=center,scale=0.8,black] at (2.5,4.5) {C};
\node [align=center,scale=0.8,black] at (3.5,4.5) {K};
\node [align=center,scale=0.8,black] at (4.5,4.5) {S};
\node [align=center,scale=0.8,black] at (0.5,1.5) {E};
\node [align=center,scale=0.8,black] at (2.5,1.5) {C};
\node [align=center,scale=0.8,black] at (4.5,1.5) {N};
\node [align=center,scale=0.8,black] at (0.5,3.5) {I};
\node [align=center,scale=0.8,black] at (2.5,3.5) {R};
\node [align=center,scale=0.8,black] at (4.5,3.5) {I};
\end{tikzpicture}
\qquad\qquad
\begin{tikzpicture}[scale=0.5] %
\draw [fill=mygreen] (0,4) rectangle (1,5);
\draw [fill=mygreen] (4,4) rectangle (5,5);
\draw [fill=mygreen] (2,0) rectangle (3,3);
\draw [fill=mygrey] (1,4) rectangle (2,5); 
\draw [fill=myyellow] (3,4) rectangle (4,5);
\draw [fill=myyellow] (2,4) rectangle (3,5); 
\draw [fill=mygrey] (0,0) rectangle (1,4);
\draw [fill=mygrey] (4,0) rectangle (5,4);
\draw [fill=mygrey] (1,0) rectangle (2,1);
\draw [fill=mygrey] (3,0) rectangle (4,1);
\draw [fill=mygrey] (1,2) rectangle (2,3);
\draw [fill=mygrey] (3,2) rectangle (4,3);
\draw [fill=mygrey] (2,1) rectangle (3,4); 
\draw [fill=mygreen] (2,2) rectangle (3,3);
\draw [fill=myyellow] (2,3) rectangle (3,4); 
\draw [-] (0,0)-- (5,0);
\draw [-] (0,1)-- (5,1);
\draw [-] (0,2)-- (5,2);
\draw [-] (0,3)-- (5,3);
\draw [-] (0,4)-- (5,4);
\draw [-] (0,5)-- (5,5);
\draw [-] (0,0)-- (0,5);
\draw [-] (1,0)-- (1,5);
\draw [-] (2,0)-- (2,5);
\draw [-] (3,0)-- (3,5);
\draw [-] (4,0)-- (4,5);
\draw [-] (5,0)-- (5,5);
\draw [fill=mywhite] (1,1) rectangle (2,2);
\draw [fill=mywhite] (1,3) rectangle (2,4);
\draw [fill=mywhite] (3,1) rectangle (4,2);
\draw [fill=mywhite] (3,3) rectangle (4,4);
\node [align=center,scale=0.8,black] at (0.5,0.5) {};
\node [align=center,scale=0.8,black] at (1.5,0.5) {};
\node [align=center,scale=0.8,black] at (2.5,0.5) {E};
\node [align=center,scale=0.8,black] at (3.5,0.5) {};
\node [align=center,scale=0.8,black] at (4.5,0.5) {};
\node [align=center,scale=0.8,black] at (0.5,2.5) {};
\node [align=center,scale=0.8,black] at (1.5,2.5) {};
\node [align=center,scale=0.8,black] at (2.5,2.5) {I};
\node [align=center,scale=0.8,black] at (3.5,2.5) {};
\node [align=center,scale=0.8,black] at (4.5,2.5) {};
\node [align=center,scale=0.8,black] at (0.5,4.5) {L};
\node [align=center,scale=0.8,black] at (1.5,4.5) {O};
\node [align=center,scale=0.8,black] at (2.5,4.5) {A};
\node [align=center,scale=0.8,black] at (3.5,4.5) {N};
\node [align=center,scale=0.8,black] at (4.5,4.5) {S};
\node [align=center,scale=0.8,black] at (0.5,1.5) {};
\node [align=center,scale=0.8,black] at (2.5,1.5) {D};
\node [align=center,scale=0.8,black] at (4.5,1.5) {};
\node [align=center,scale=0.8,black] at (0.5,3.5) {};
\node [align=center,scale=0.8,black] at (2.5,3.5) {S};
\node [align=center,scale=0.8,black] at (4.5,3.5) {};
\end{tikzpicture}
\vskip-2mm
\caption{(a) (S,T) gives BOOTS, NASAL, (b) (A,I,O) gives
  BISON, SACKS, (c) (O,U,Y) gives SUCKS, DOING, DYKES and (d) (O,S,D,N,A)
  gives LANDS, NOISE}
\label{F:2}
\end{center}
\end{figure}

An implicit (and important) rule of Waffle is that there is \emph{only} one answer, and hence
only one perfect unscrambling. (If there is a unique unscrambling $g\in S_{21}$,
there are still many different sequences of 10 swaps giving $g$. For example,
an $n_i$-cycle with $n_i\ge2$ is a product of $n_i-1$
swaps in $n_i$ different ways.)
Call a 5-letter English word \emph{swappable} if swapping the two even
squares gives a different 5-letter English word. For example, applying
the swap $({\rm M},{\rm L})$ to
\begin{tikzpicture}[scale=0.4] %
\draw [fill=mygreen] (0,0) rectangle (1,1);
\draw [fill=mygreen] (2,0) rectangle (3,1);
\draw [fill=mygreen] (4,0) rectangle (5,1);
\draw [fill=myyellow] (1,0) rectangle (2,1);
\draw [fill=myyellow] (3,0) rectangle (4,1);
\draw [-] (0,0) -- (5,0) -- (5,1) -- (0,1) -- (0,0);
\draw [-] (1,0) -- (1,1);
\draw [-] (2,0) -- (2,1);
\draw [-] (3,0) -- (3,1);
\draw [-] (4,0) -- (4,1);
\node [align=center,scale=0.8,black] at (0.5,0.5) {S};
\node [align=center,scale=0.8,black] at (1.5,0.5) {M};
\node [align=center,scale=0.8,black] at (2.5,0.5) {I};
\node [align=center,scale=0.8,black] at (3.5,0.5) {L};
\node [align=center,scale=0.8,black] at (4.5,0.5) {E};
\end{tikzpicture}
gives the different English word 
\begin{tikzpicture}[scale=0.4] %
\draw [fill=mygreen] (0,0) rectangle (5,1);
\draw [-] (0,0) -- (5,0) -- (5,1) -- (0,1) -- (0,0);
\draw [-] (1,0) -- (1,1);
\draw [-] (2,0) -- (2,1);
\draw [-] (3,0) -- (3,1);
\draw [-] (4,0) -- (4,1);
\node [align=center,scale=0.8,black] at (0.5,0.5) {S};
\node [align=center,scale=0.8,black] at (1.5,0.5) {L};
\node [align=center,scale=0.8,black] at (2.5,0.5) {I};
\node [align=center,scale=0.8,black] at (3.5,0.5) {M};
\node [align=center,scale=0.8,black] at (4.5,0.5) {E};
\end{tikzpicture}.
We exclude swappable words from our dictionary of possible 5-letter words.
Similarly, we must exclude certain pairs and
triples of words, etc. such as in Figure~\ref{F:2}.
The archived~\cite{WA} Waffle game \#679 sadly does not have a unique solution:
the swap $({\rm S},{\rm T})$ shown in Figure~\ref{F:2}(a) changes one valid
solution to another. (Both can be solved in 10 swaps.)

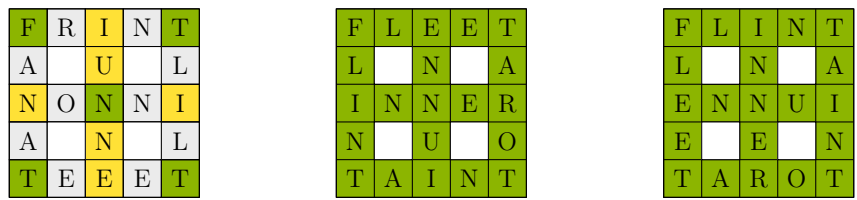
\begin{figure}[!ht]
\begin{center}
\begin{tikzpicture}[scale=0.5] %
\draw [fill=mygreen] (0,0) rectangle (1,1);
\draw [fill=mygreen] (4,0) rectangle (5,1);
\draw [fill=mygreen] (0,4) rectangle (1,5);
\draw [fill=mygreen] (4,4) rectangle (5,5);
\draw [fill=mygreen] (2,2) rectangle (3,3);
\draw [fill=myyellow] (2,0) rectangle (3,1);
\draw [fill=myyellow] (4,1) rectangle (5,4);
\draw [fill=myyellow] (1,2) rectangle (2,3);
\draw [fill=myyellow] (0,3) rectangle (1,4);
\draw [fill=myyellow] (2,4) rectangle (3,5);
\draw [fill=mygrey] (1,0) rectangle (2,1);
\draw [fill=mygrey] (3,0) rectangle (4,1);
\draw [fill=mygrey] (0,1) rectangle (1,3);
\draw [fill=mygrey] (2,1) rectangle (3,2);
\draw [fill=mygrey] (3,2) rectangle (4,3);
\draw [fill=mygrey] (2,3) rectangle (3,4);
\draw [fill=mygrey] (1,4) rectangle (2,5);
\draw [fill=mygrey] (3,4) rectangle (4,5);
\draw [-] (0,0)-- (5,0);
\draw [-] (0,1)-- (5,1);
\draw [-] (0,2)-- (5,2);
\draw [-] (0,3)-- (5,3);
\draw [-] (0,4)-- (5,4);
\draw [-] (0,5)-- (5,5);
\draw [-] (0,0)-- (0,5);
\draw [-] (1,0)-- (1,5);
\draw [-] (2,0)-- (2,5);
\draw [-] (3,0)-- (3,5);
\draw [-] (4,0)-- (4,5);
\draw [-] (5,0)-- (5,5);
\draw [fill=mywhite] (1,1) rectangle (2,2);
\draw [fill=mywhite] (1,3) rectangle (2,4);
\draw [fill=mywhite] (3,1) rectangle (4,2);
\draw [fill=mywhite] (3,3) rectangle (4,4);
\node [align=center,scale=0.8,black] at (0.5,0.5) {N};
\node [align=center,scale=0.8,black] at (1.5,0.5) {N};
\node [align=center,scale=0.8,black] at (2.5,0.5) {B};
\node [align=center,scale=0.8,black] at (3.5,0.5) {S};
\node [align=center,scale=0.8,black] at (4.5,0.5) {L};
\node [align=center,scale=0.8,black] at (0.5,2.5) {O};
\node [align=center,scale=0.8,black] at (1.5,2.5) {A};
\node [align=center,scale=0.8,black] at (2.5,2.5) {O};
\node [align=center,scale=0.8,black] at (3.5,2.5) {I};
\node [align=center,scale=0.8,black] at (4.5,2.5) {O};
\node [align=center,scale=0.8,black] at (0.5,4.5) {C};
\node [align=center,scale=0.8,black] at (1.5,4.5) {G};
\node [align=center,scale=0.8,black] at (2.5,4.5) {A};
\node [align=center,scale=0.8,black] at (3.5,4.5) {T};
\node [align=center,scale=0.8,black] at (4.5,4.5) {L};
\node [align=center,scale=0.8,black] at (0.5,1.5) {A};
\node [align=center,scale=0.8,black] at (2.5,1.5) {A};
\node [align=center,scale=0.8,black] at (4.5,1.5) {Y};
\node [align=center,scale=0.8,black] at (0.5,3.5) {A};
\node [align=center,scale=0.8,black] at (2.5,3.5) {H};
\node [align=center,scale=0.8,black] at (4.5,3.5) {A};
\end{tikzpicture}
\qquad\qquad
\begin{tikzpicture}[scale=0.5] %
\draw [fill=mygreen] (0,0) rectangle (5,5);
\draw [-] (0,0)-- (5,0);
\draw [-] (0,1)-- (5,1);
\draw [-] (0,2)-- (5,2);
\draw [-] (0,3)-- (5,3);
\draw [-] (0,4)-- (5,4);
\draw [-] (0,5)-- (5,5);
\draw [-] (0,0)-- (0,5);
\draw [-] (1,0)-- (1,5);
\draw [-] (2,0)-- (2,5);
\draw [-] (3,0)-- (3,5);
\draw [-] (4,0)-- (4,5);
\draw [-] (5,0)-- (5,5);
\draw [fill=mywhite] (1,1) rectangle (2,2);
\draw [fill=mywhite] (1,3) rectangle (2,4);
\draw [fill=mywhite] (3,1) rectangle (4,2);
\draw [fill=mywhite] (3,3) rectangle (4,4);
\node [align=center,scale=0.8,black] at (0.5,0.5) {N};
\node [align=center,scale=0.8,black] at (1.5,0.5) {A};
\node [align=center,scale=0.8,black] at (2.5,0.5) {T};
\node [align=center,scale=0.8,black] at (3.5,0.5) {A};
\node [align=center,scale=0.8,black] at (4.5,0.5) {L};
\node [align=center,scale=0.8,black] at (0.5,2.5) {A};
\node [align=center,scale=0.8,black] at (1.5,2.5) {G};
\node [align=center,scale=0.8,black] at (2.5,2.5) {O};
\node [align=center,scale=0.8,black] at (3.5,2.5) {N};
\node [align=center,scale=0.8,black] at (4.5,2.5) {Y};
\node [align=center,scale=0.8,black] at (0.5,4.5) {C};
\node [align=center,scale=0.8,black] at (1.5,4.5) {A};
\node [align=center,scale=0.8,black] at (2.5,4.5) {B};
\node [align=center,scale=0.8,black] at (3.5,4.5) {A};
\node [align=center,scale=0.8,black] at (4.5,4.5) {L};
\node [align=center,scale=0.8,black] at (0.5,1.5) {I};
\node [align=center,scale=0.8,black] at (2.5,1.5) {S};
\node [align=center,scale=0.8,black] at (4.5,1.5) {A};
\node [align=center,scale=0.8,black] at (0.5,3.5) {H};
\node [align=center,scale=0.8,black] at (2.5,3.5) {O};
\node [align=center,scale=0.8,black] at (4.5,3.5) {O};
\end{tikzpicture}
\qquad\qquad
\begin{tikzpicture}[scale=0.5] %
\draw [fill=mygreen] (0,0) rectangle (5,5);
\draw [fill=myyellow] (2,0) rectangle (3,2);
\draw [-] (0,0)-- (5,0);
\draw [-] (0,1)-- (5,1);
\draw [-] (0,2)-- (5,2);
\draw [-] (0,3)-- (5,3);
\draw [-] (0,4)-- (5,4);
\draw [-] (0,5)-- (5,5);
\draw [-] (0,0)-- (0,5);
\draw [-] (1,0)-- (1,5);
\draw [-] (2,0)-- (2,5);
\draw [-] (3,0)-- (3,5);
\draw [-] (4,0)-- (4,5);
\draw [-] (5,0)-- (5,5);
\draw [fill=mywhite] (1,1) rectangle (2,2);
\draw [fill=mywhite] (1,3) rectangle (2,4);
\draw [fill=mywhite] (3,1) rectangle (4,2);
\draw [fill=mywhite] (3,3) rectangle (4,4);
\node [align=center,scale=0.8,black] at (0.5,0.5) {N};
\node [align=center,scale=0.8,black] at (1.5,0.5) {A};
\node [align=center,scale=0.8,black] at (2.5,0.5) {S};
\node [align=center,scale=0.8,black] at (3.5,0.5) {A};
\node [align=center,scale=0.8,black] at (4.5,0.5) {L};
\node [align=center,scale=0.8,black] at (0.5,2.5) {A};
\node [align=center,scale=0.8,black] at (1.5,2.5) {G};
\node [align=center,scale=0.8,black] at (2.5,2.5) {O};
\node [align=center,scale=0.8,black] at (3.5,2.5) {N};
\node [align=center,scale=0.8,black] at (4.5,2.5) {Y};
\node [align=center,scale=0.8,black] at (0.5,4.5) {C};
\node [align=center,scale=0.8,black] at (1.5,4.5) {A};
\node [align=center,scale=0.8,black] at (2.5,4.5) {B};
\node [align=center,scale=0.8,black] at (3.5,4.5) {A};
\node [align=center,scale=0.8,black] at (4.5,4.5) {L};
\node [align=center,scale=0.8,black] at (0.5,1.5) {I};
\node [align=center,scale=0.8,black] at (2.5,1.5) {T};
\node [align=center,scale=0.8,black] at (4.5,1.5) {A};
\node [align=center,scale=0.8,black] at (0.5,3.5) {H};
\node [align=center,scale=0.8,black] at (2.5,3.5) {O};
\node [align=center,scale=0.8,black] at (4.5,3.5) {O};
\end{tikzpicture}
\vskip-2mm
\caption{(a) archived Waffle game \# 679 sadly has two valid solutions~(b,c)}
\label{F:2.5}
\end{center}
\end{figure}

We now consider algorithms which solve Step~1, and determine whether
a Waffle game has a unique solution.
One completely impractical algorithm involves checking whether each
non-identity $g\in S_{21}$ applied to the set of six (scrambled) 5-letter
words gives a set of six (unscrambled) English words. This is impractical
as there are $21!-1\approx 10^{10}$ choices for $g$.
Another algorithm which is practical, finds for each  horizontal word
$h_1,h_2,h_3$ and each vertical word $v_1,v_2,v_3$ a list of candidate
words from our dictionary of words that satisfy the color constraints; and then
determines which 6 tuples $(h_1,h_2,h_3,v_1,v_2,v_3)$ have the same letters,
with the same multiplicities, as the scrambled game. This algorithm
is implemented in~\cite{G} using {\sc Magma}~\cite{M} code.
It is suited for a computer and not
a person. Using this program we found that Waffle Game \#~679 in
Figure~\ref{F:2.5}(a) had the two solutions in Figure~\ref{F:2.5}(b,c).
This program works well on the archived Waffle games~\cite{WA} with at
least 5~green squares, but is struggles with the Waffle game in
Figure~\ref{F:3}(c) with only one green square.

\section{What makes some games hard?}\label{S:hard}

As you play more games of Waffle, you will notice that some games are
easier to achieve a perfect score than others.  This section considers
what features of a Waffle game makes it easy or hard. Players want
`Goldilocks games', which are `not too easy, not too hard, but just right'.
Thus design features must be made such as choosing a list of 5-letter
English words (e.g.~\cite{W}) that is large enough, but not so large
as to include unfamiliar words. I initially used the valid 5-letter
Scrabble words, but this list has many unfamiliar words (e.g. ABUNA, BLAER, COMTE, EEJIT, etc).

\begin{figure}[!ht]
\begin{center}
\begin{tikzpicture}[scale=0.5] %
\draw [fill=mygrey] (0,0) rectangle (5,5);
\draw [fill=myyellow] (0,0) rectangle (1,4);
\draw [fill=myyellow] (1,4) rectangle (5,5);
\draw [fill=myyellow] (1,0) rectangle (5,1);
\draw [fill=myyellow] (4,1) rectangle (5,5);
\draw [fill=myyellow] (0,2) rectangle (5,3);
\draw [fill=myyellow] (2,1) rectangle (3,2);
\draw [fill=myyellow] (2,3) rectangle (3,4);
\draw [fill=mygreen] (0,4) rectangle (1,5);
\draw [fill=mywhite] (1,1) rectangle (2,2);
\draw [fill=mywhite] (1,3) rectangle (2,4);
\draw [fill=mywhite] (3,1) rectangle (4,2);
\draw [fill=mywhite] (3,3) rectangle (4,4);
\draw [fill=mygrey] (2,2) rectangle (3,3);
\draw [fill=mygrey] (3,0) rectangle (5,1);
\draw [-] (0,0)-- (5,0);
\draw [-] (0,1)-- (5,1);
\draw [-] (0,2)-- (5,2);
\draw [-] (0,3)-- (5,3);
\draw [-] (0,4)-- (5,4);
\draw [-] (0,5)-- (5,5);
\draw [-] (0,0)-- (0,5);
\draw [-] (1,0)-- (1,5);
\draw [-] (2,0)-- (2,5);
\draw [-] (3,0)-- (3,5);
\draw [-] (4,0)-- (4,5);
\draw [-] (5,0)-- (5,5);
\node [align=center,scale=0.8,black] at (0.5,0.5) {M};
\node [align=center,scale=0.8,black] at (1.5,0.5) {S};
\node [align=center,scale=0.8,black] at (2.5,0.5) {M};
\node [align=center,scale=0.8,black] at (3.5,0.5) {M};
\node [align=center,scale=0.8,black] at (4.5,0.5) {M};
\node [align=center,scale=0.8,black] at (0.5,2.5) {A};
\node [align=center,scale=0.8,black] at (1.5,2.5) {M};
\node [align=center,scale=0.8,black] at (2.5,2.5) {A};
\node [align=center,scale=0.8,black] at (3.5,2.5) {A};
\node [align=center,scale=0.8,black] at (4.5,2.5) {S};
\node [align=center,scale=0.8,black] at (0.5,4.5) {M};
\node [align=center,scale=0.8,black] at (1.5,4.5) {M};
\node [align=center,scale=0.8,black] at (2.5,4.5) {A};
\node [align=center,scale=0.8,black] at (3.5,4.5) {A};
\node [align=center,scale=0.8,black] at (4.5,4.5) {M};
\node [align=center,scale=0.8,black] at (0.5,1.5) {A};
\node [align=center,scale=0.8,black] at (2.5,1.5) {A};
\node [align=center,scale=0.8,black] at (4.5,1.5) {A};
\node [align=center,scale=0.8,black] at (0.5,3.5) {M};
\node [align=center,scale=0.8,black] at (2.5,3.5) {M};
\node [align=center,scale=0.8,black] at (4.5,3.5) {S};
\end{tikzpicture}
\quad
\begin{tikzpicture}[scale=0.5] %
\draw [fill=mygrey] (0,0) rectangle (5,5);
\draw [fill=mygreen] (0,0) rectangle (5,1);
\draw [fill=mygreen] (0,2) rectangle (5,3);
\draw [fill=mygreen] (0,4) rectangle (5,5);
\draw [fill=mygreen] (0,0) rectangle (1,5);
\draw [fill=mygreen] (2,0) rectangle (3,5);
\draw [fill=mygreen] (4,0) rectangle (5,5);
\draw [fill=mywhite] (1,1) rectangle (2,2);
\draw [fill=mywhite] (1,3) rectangle (2,4);
\draw [fill=mywhite] (3,1) rectangle (4,2);
\draw [fill=mywhite] (3,3) rectangle (4,4);
\draw [-] (0,0)-- (5,0);
\draw [-] (0,1)-- (5,1);
\draw [-] (0,2)-- (5,2);
\draw [-] (0,3)-- (5,3);
\draw [-] (0,4)-- (5,4);
\draw [-] (0,5)-- (5,5);
\draw [-] (0,0)-- (0,5);
\draw [-] (1,0)-- (1,5);
\draw [-] (2,0)-- (2,5);
\draw [-] (3,0)-- (3,5);
\draw [-] (4,0)-- (4,5);
\draw [-] (5,0)-- (5,5);
\node [align=center,scale=0.8,black] at (0.5,0.5) {A};
\node [align=center,scale=0.8,black] at (1.5,0.5) {M};
\node [align=center,scale=0.8,black] at (2.5,0.5) {A};
\node [align=center,scale=0.8,black] at (3.5,0.5) {S};
\node [align=center,scale=0.8,black] at (4.5,0.5) {S};
\node [align=center,scale=0.8,black] at (0.5,2.5) {M};
\node [align=center,scale=0.8,black] at (1.5,2.5) {A};
\node [align=center,scale=0.8,black] at (2.5,2.5) {M};
\node [align=center,scale=0.8,black] at (3.5,2.5) {M};
\node [align=center,scale=0.8,black] at (4.5,2.5) {A};
\node [align=center,scale=0.8,black] at (0.5,4.5) {M};
\node [align=center,scale=0.8,black] at (1.5,4.5) {A};
\node [align=center,scale=0.8,black] at (2.5,4.5) {M};
\node [align=center,scale=0.8,black] at (3.5,4.5) {M};
\node [align=center,scale=0.8,black] at (4.5,4.5) {A};
\node [align=center,scale=0.8,black] at (0.5,1.5) {M};
\node [align=center,scale=0.8,black] at (2.5,1.5) {M};
\node [align=center,scale=0.8,black] at (4.5,1.5) {S};
\node [align=center,scale=0.8,black] at (0.5,3.5) {A};
\node [align=center,scale=0.8,black] at (2.5,3.5) {A};
\node [align=center,scale=0.8,black] at (4.5,3.5) {M};
\end{tikzpicture}
\quad\quad\quad
\begin{tikzpicture}[scale=0.5] %
\draw [fill=mygrey] (0,0) rectangle (5,5);
\draw [-] (0,0)-- (5,0);
\draw [-] (0,1)-- (5,1);
\draw [-] (0,2)-- (5,2);
\draw [-] (0,3)-- (5,3);
\draw [-] (0,4)-- (5,4);
\draw [-] (0,5)-- (5,5);
\draw [-] (0,0)-- (0,5);
\draw [-] (1,0)-- (1,5);
\draw [-] (2,0)-- (2,5);
\draw [-] (3,0)-- (3,5);
\draw [-] (4,0)-- (4,5);
\draw [-] (5,0)-- (5,5);
\draw [fill=mygreen] (2,3) rectangle (3,4);
\draw [fill=mywhite] (1,1) rectangle (2,2);
\draw [fill=mywhite] (1,3) rectangle (2,4);
\draw [fill=mywhite] (3,1) rectangle (4,2);
\draw [fill=mywhite] (3,3) rectangle (4,4);
\node [align=center,scale=0.8,black] at (0.5,0.5) {T};
\node [align=center,scale=0.8,black] at (1.5,0.5) {L};
\node [align=center,scale=0.8,black] at (2.5,0.5) {H};
\node [align=center,scale=0.8,black] at (3.5,0.5) {M};
\node [align=center,scale=0.8,black] at (4.5,0.5) {B};
\node [align=center,scale=0.8,black] at (0.5,2.5) {E};
\node [align=center,scale=0.8,black] at (1.5,2.5) {W};
\node [align=center,scale=0.8,black] at (2.5,2.5) {K};
\node [align=center,scale=0.8,black] at (3.5,2.5) {R};
\node [align=center,scale=0.8,black] at (4.5,2.5) {G};
\node [align=center,scale=0.8,black] at (0.5,4.5) {D};
\node [align=center,scale=0.8,black] at (1.5,4.5) {U};
\node [align=center,scale=0.8,black] at (2.5,4.5) {P};
\node [align=center,scale=0.8,black] at (3.5,4.5) {O};
\node [align=center,scale=0.8,black] at (4.5,4.5) {S};
\node [align=center,scale=0.8,black] at (0.5,1.5) {A};
\node [align=center,scale=0.8,black] at (2.5,1.5) {Y};
\node [align=center,scale=0.8,black] at (4.5,1.5) {C};
\node [align=center,scale=0.8,black] at (0.5,3.5) {I};
\node [align=center,scale=0.8,black] at (2.5,3.5) {U$'$};
\node [align=center,scale=0.8,black] at (4.5,3.5) {N};
\end{tikzpicture}
\quad
\begin{tikzpicture}[scale=0.5] %
\draw [fill=mygrey] (0,0) rectangle (5,5);
\draw [fill=mygreen] (0,0) rectangle (5,1);
\draw [fill=mygreen] (0,2) rectangle (5,3);
\draw [fill=mygreen] (0,4) rectangle (5,5);
\draw [fill=mygreen] (0,0) rectangle (1,5);
\draw [fill=mygreen] (2,0) rectangle (3,5);
\draw [fill=mygreen] (4,0) rectangle (5,5);
\draw [-] (0,0)-- (5,0);
\draw [-] (0,1)-- (5,1);
\draw [-] (0,2)-- (5,2);
\draw [-] (0,3)-- (5,3);
\draw [-] (0,4)-- (5,4);
\draw [-] (0,5)-- (5,5);
\draw [-] (0,0)-- (0,5);
\draw [-] (1,0)-- (1,5);
\draw [-] (2,0)-- (2,5);
\draw [-] (3,0)-- (3,5);
\draw [-] (4,0)-- (4,5);
\draw [-] (5,0)-- (5,5);
\draw [fill=mywhite] (1,1) rectangle (2,2);
\draw [fill=mywhite] (1,3) rectangle (2,4);
\draw [fill=mywhite] (3,1) rectangle (4,2);
\draw [fill=mywhite] (3,3) rectangle (4,4);
\node [align=center,scale=0.8,black] at (0.5,0.5) {S};
\node [align=center,scale=0.8,black] at (1.5,0.5) {W};
\node [align=center,scale=0.8,black] at (2.5,0.5) {O};
\node [align=center,scale=0.8,black] at (3.5,0.5) {R};
\node [align=center,scale=0.8,black] at (4.5,0.5) {D};
\node [align=center,scale=0.8,black] at (0.5,2.5) {C};
\node [align=center,scale=0.8,black] at (1.5,2.5) {L};
\node [align=center,scale=0.8,black] at (2.5,2.5) {A};
\node [align=center,scale=0.8,black] at (3.5,2.5) {M};
\node [align=center,scale=0.8,black] at (4.5,2.5) {P};
\node [align=center,scale=0.8,black] at (0.5,4.5) {B};
\node [align=center,scale=0.8,black] at (1.5,4.5) {I};
\node [align=center,scale=0.8,black] at (2.5,4.5) {G};
\node [align=center,scale=0.8,black] at (3.5,4.5) {H};
\node [align=center,scale=0.8,black] at (4.5,4.5) {T};
\node [align=center,scale=0.8,black] at (0.5,1.5) {K};
\node [align=center,scale=0.8,black] at (2.5,1.5) {N};
\node [align=center,scale=0.8,black] at (4.5,1.5) {E};
\node [align=center,scale=0.8,black] at (0.5,3.5) {U};
\node [align=center,scale=0.8,black] at (2.5,3.5) {U$'$};
\node [align=center,scale=0.8,black] at (4.5,3.5) {Y};
\end{tikzpicture}
\end{center}
\vskip-4mm
\caption{$N_g=1$ (a) easy scrambled, (b) unscrambled, (c) hard, (d) unscrambled}
\label{F:3}
\end{figure}
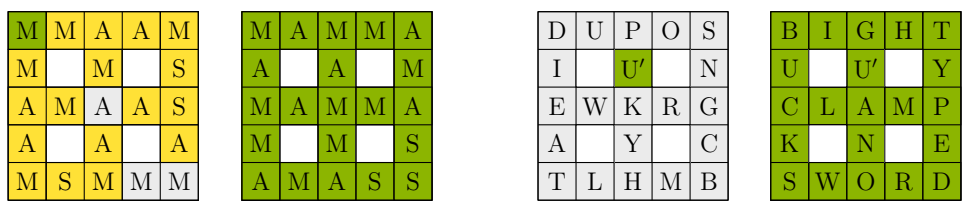

We saw in Section~\ref{S:color} that many repeated letters can make it
harder to find a perfect unscrambling via a naive method. (However,
heuristics using the fact that a perfect unscrambling has  a large number of
cycles can vastly reduce the possibilities.)
Finding 10 swaps from the
perfect unscrambling (Step~2) is then easy by the proof of Cayley's lemma, but finding a perfect unscrambling (Step~1)
may be hard without a computer program. Although games with many repeated
letters can require heuristics to reduce the space of possibilities, this is
not always true: the scrambled
game Figure~\ref{F:3}(a) has 10 Ms, 8 As, and 3 Ss, and yet it is
easy to solve. This Waffle game is unusual in three ways:
(a) it has only one green square, (b) Figure~\ref{F:3}(b) has \emph{two}
distinct words and most Waffle games have six distinct words, and (c) of the
$\binom{21}{10}\binom{11}{8}\binom{3}{3}=\frac{21!}{10!\,8!\,3!}\approx 10^{8}$ choices for $g$ many give perfect unscramblings.

The difficulty of a given Waffle game (before any swaps)
seems to depend on the numbers $N_g, N_y$ of green and yellow
squares, respectively. The number of gray square is $21-N_g-N_y$.  We proved
in Section~\ref{S:color} that $1\le N_g\le10$. The goal of Waffle is to increase
$N_g$ to 21 with 10 swaps. Clearly if $N_g=19$, then the last swap
will give $N_g=21$, but sometimes the best that you can do with one swap is to
increase $N_g$ by one. Hence our intuition is that finding a perfect
unscrambling for a puzzle with $N_g=1$ should be particularly hard,
especially if you are not using a computer for this task.

Waffle games that give you a lot of information about the mystery six words
have larger values of $N_g$, and to a lesser extent $N_y$. Gray squares
give information too, but less. Hence a
Waffle game with $(N_g,N_y)=(1,0)$ seems to give as little
information as possible. Let $g$ be the perfect unscrambling for such a
game, and let $g'$ be the permutation of the 20 non-green squares. Since
$c(g')=c(g)-1=10$, the cycle decomposition of $g'$ must
be $g'=g'_1\cdots g'_{10}$ where $g'_1,\dots,g'_{10}$ are \emph{disjoint}
2-cycles (i.e. a swaps).
The 2-cycles $g'_1,\dots,g'_{10}$ can be permuted in $10!$ different ways, and we
still obtain a valid disjoint cycle decomposition for $g'$.
Assume that $g'$ permutes 20 \emph{distinct} letters (without subscripts) as in
Figure~\ref{F:3}(c). This implies that only one choice for
$g'$ gives a perfect unscrambling, and $c(g')=10$ must hold. The number of
permutations $g'\in S_{20}$ with a disjoint cycle decomposition comprising
ten 2-cycles is
$\frac{20!}{(2!)^{10}10!}\approx 10^{8.8}$, so the chance of finding the right
$g'$ at random is $\approx 10^{-8.8}$. 
This seems impossibly hard without a computer program. A naive program
running through all $\approx 10^{8.8}$ possibilities would certainly not terminate. However,
a smarter program using statistical properties of 5-letter English words may
be able to reduce the search space from $\approx 10^{8.8}$ to a manageable size.
Of the $26^2=676$ two letter strings $xy$ with
$x,y\in\{\textup{A,\dots,Z}\}$ precisely 306 choices for $xy$ have
frequency zero in list of 3075 common non-swappable 5-letter words
found from the list~\cite{W}.

In summary, a Waffle game such as Figure~\ref{F:3}(c)
with $(N_g,N_y)=(1,0)$ and 20 distinct non-green letters
is extremely hard to solve by guessing,
and may even be too hard for a clever computer program.
This gives some insight why the the creator of Waffle,
James Robinson, creates games with $N_g\ge5$: these games are easier for humans
to solve. For the Waffle game in Figure~\ref{F:3}(c) finding the six 5-letter
words (Step~1) in Figure~\ref{F:3}(d) is the hard part. Finding a
perfect unscrambling (Step~2)
is then easy: there is only one choice namely
$g=\textup{(B,D)(I,U)(G,P)(H,O)(T,S)(Y,N)(C,E)(L,W)(A,K)(M,R)(U$'$)}$.
Finally, we have already found the 10 swaps so Step~3 is vacuous.

\begin{thebibliography}{9}

\bibitem{M}
  Wieb Bosma, John Cannon and Catherine Playoust,
  The Magma algebra system. I. The user language.
  Computational algebra and number theory (London, 1993).
  \emph{J. Symbolic Comput.} {\bf24} (1997), no. 3-4, 235--265.

\bibitem{JR}  Christopher Dring, I've turned down life-changing money: The story of Waffle,\newline
  \url{https://tinyurl.com/bvsxs9h3}

\bibitem{ET2}
  P. Erd\H{o}s and P. Tur\'{a}n,
  On some problems of a statistical group-theory. I,
  \emph{Z. Wahrscheinlichkeitstheorie und Verw. Gebiete},
   (1965){\bf4}, 175--186.

\bibitem{G}
  S.P. Glasby.
  {\sc Magma} computer code for computing with Waffle,\newline
  \url{https://stephenglasby.github.io/Waffle/}


\bibitem{CM} Carl Mueller, Personal communication, 30 January 2025.

\bibitem{WA}
  Waffle archive, \url{https://wafflegame.net/archive}

\bibitem{W}
  Luke Williams, 3103 common 5-letter words,
  \url{https://tinyurl.com/46rhrc5d}

\end{thebibliography}
\end{document}